\documentclass[12pt,leqno,twoside]{amsart}

\input boxedeps   
\SetRokickiEPSFSpecial 
\HideDisplacementBoxes                  
\newcommand{\picspace}{\vspace{4mm}}    

\makeatletter
\def\thmhead@plain#1#2#3{\thmname{#1}}

\newtheorem{theorem}{Theorem} 
\newtheorem{corollary}[theorem]{Corollary}

\def\cmc/{\textsc{cmc}}
\newcommand{\R}{\mathbb R}
\renewcommand{\S}{\mathbb S}

\renewcommand{\baselinestretch}{1.12} 

\begin{document}

\author[Grosse-Brauckmann]{Karsten Gro\ss e-Brauckmann}
\address{Universit\"at Bonn\\ Mathematisches Institut\\ Beringstr.\ 1\\
  53115 Bonn, Germany}
  \email{kgb@math.uni-bonn.de}

\author[Kusner]{Robert B. Kusner}
\address{Mathematics Department, University of Massachusetts,
  Amherst MA 01003, USA}
  \email{kusner@math.umass.edu}                                 

\author[Sullivan]{John M. Sullivan}
\address{Mathematics Department, University of Illinois,
  Urbana IL 61801, USA}
  \email{sullivan@math.uiuc.edu}

\title[Constant mean curvature surfaces with three ends]
{Constant mean curvature surfaces\\ with three ends}

\begin{abstract}
We announce the classification of complete,
almost embedded surfaces of constant mean curvature, with three ends
and genus zero: they are classified by triples of points on the sphere
whose distances are the asymptotic necksizes of the three ends.
%
\end{abstract}

\maketitle \thispagestyle{empty}


Surfaces which minimize area under a volume constraint have constant mean
curvature (\cmc/); this condition can be expressed as a nonlinear partial
differential equation.  We are interested in complete \cmc/ surfaces
properly embedded in $\R^3$; we rescale them to have mean curvature one.
For technical reasons, we consider a slight generalization of embeddedness
(introduced by Alexandrov~\cite{al2}): An immersed surface is \emph{almost
embedded} if it bounds a properly immersed three-manifold.

Alexandrov~\cite{al1,al2} showed that the round sphere is the only
\emph{compact} almost embedded \cmc/ surface.  The next case to consider is
that of \emph{finite-topology} surfaces, homeomorphic to a compact surface
with a finite number of points removed.  A neighborhood of any of these
punctures is called an \emph{end} of the surface.  The \emph{unduloids},
\cmc/ surfaces of revolution described by Delaunay~\cite{del}, are
genus-zero examples with two ends.  Each is a solution of an ordinary
differential equation; the entire family is parametrized by the unduloid
\emph{necksize}, which ranges from zero (at the singular chain of spheres)
to $\pi$ (at the cylinder).

Over the past decade there has been increasing understanding of
finite-topology almost embedded \cmc/ surfaces.  Each end of such a surface
is asymptotic to an unduloid~\cite{kks}.  Meeks~\cite{mee} showed there
are no examples with a single end.  The unduloids themselves are the only
examples with two ends~\cite{kks}.  Kapouleas~\cite{kap} has constructed
examples (near the zero-necksize limit) with any genus and any number of
ends greater than two.

In this note we announce the classification of all almost embedded \cmc/
surfaces with three ends and genus zero; we call these \emph{triunduloids}
(see figure).  In light of the trousers decomposition for surfaces,
triunduloids can be seen as the building blocks for more complicated
almost embedded \cmc/ surfaces~\cite{gbks2}.  Our main result determines
explicitly the \emph{moduli space} of triunduloids with labelled ends, up
to Euclidean motions.  Since triunduloids are transcendental objects, and
are not described by any ordinary differential equation, it is remarkable
to have such a complete and explicit determination for their moduli space.

\begin{figure}[th]
\ForceWidth{10cm} \centerline{\BoxedEPSF{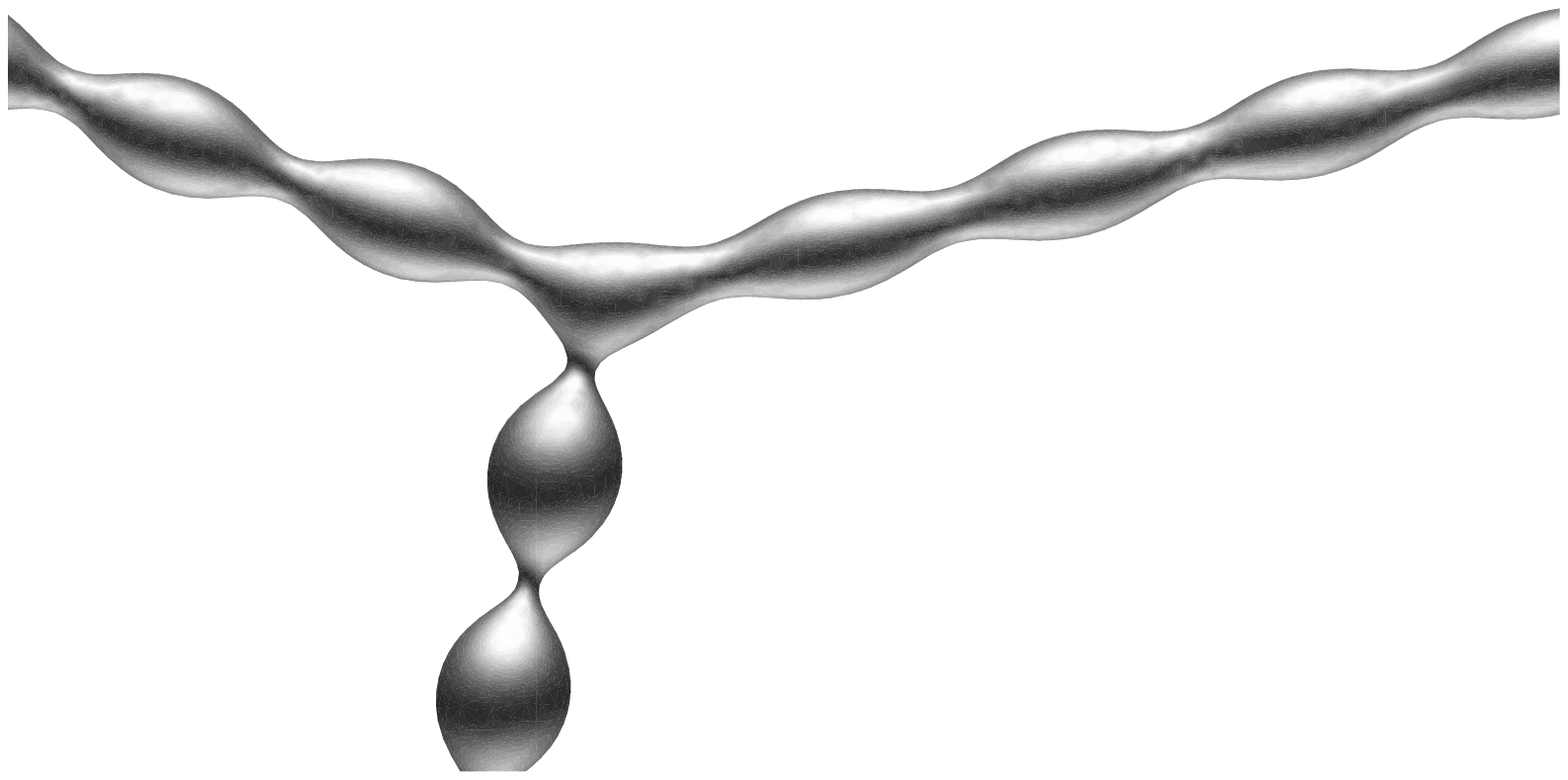}} \picspace
\centerline{\textbf{Figure.}
A triunduloid is an embedded constant mean curvature surface}
\centerline{with three ends, each asymptotic to a Delaunay unduloid.}
\end{figure}

\begin{theorem}
Triunduloids are classified by triples of distinct labelled points in the
two-sphere (up to rotations); the spherical distances of points in the
triple are the necksizes of the unduloids asymptotic to the three ends.
The moduli space of triunduloids is therefore homeomorphic to an open
three-ball.
\end{theorem}

The proof of the theorem has three parts.  First we define the classifying
map from triunduloids to spherical triples, and observe that it is proper;
then we prove it is injective; and finally we show it is surjective.

To define the classifying map, we use the fact that any triunduloid has
a reflection symmetry which decomposes the surface into mirror-image
halves~\cite{gk}.  Each half is simply connected, so Lawson's
construction~\cite{law} gives a conjugate cousin minimal surface in the
three-sphere.  Using observations of Karcher~\cite{kar} we find that its
boundary projects under the Hopf map to the desired spherical triple.
The composition of these steps defines our classifying map~\cite{gbks1}.
It follows from curvature estimates~\cite{kk} that the map is proper.

The injectivity of our classifying map is really a uniqueness result.
We use the Hopf circle bundle to construct a trivial circle bundle over
the disk representing the Lawson conjugate.  Its total space is locally
isometric to the three-sphere, and so the circle action along the fibres
is by isometries.  The classifying triple determines the bundle up to
isometries.  Moreover, any conjugate surface with the same triple defines
a minimal section of the bundle.  Thus we are in a situation familiar
from minimal graphs, and we can apply a suitable maximum principle to
deduce uniqueness.

Finally, we need an existence result showing that our classifying map is
surjective.  We depend on the fact~\cite{kmp} that the moduli space of \cmc/
surfaces of genus~$g$ with $k$~ends is locally a real analytic variety of
(formal) dimension $3k-6$.   In particular, near a \emph{nondegenerate}
triunduloid, our moduli space has dimension three.

The intersection theory for real analytic varieties developed by Borel
and Haefliger~\cite{bh} allows us to define the mapping degree of our
classifying map.  It counts (mod~2) the number of points in any generic
fiber of the map (meaning a fiber consisting entirely of nondegenerate
triunduloids).  By our uniqueness result, each fiber contains at most one
triunduloid.  This means that the existence of a single nondegenerate
triunduloid guarantees that the degree (mod~2) is~1, and thus that
every fiber contains a triunduloid.  Using a nondegenerate minimal
surface of Montiel and Ros~\cite{mr} in a construction by Mazzeo and
Pacard~\cite{mp}, we can obtain a triunduloid known to be nondegenerate.
This proves surjectivity.

In the category of real analytic varieties, a proper bijective map is
actually a homeomorphism, so the proof of the theorem is complete.
\qed\bigskip

Note that our geometric picture of the triunduloid moduli space naturally
explains necksize bounds for triunduloids.  For instance, the symmetric
triunduloids constructed previously~\cite{kgb} have three congruent ends,
with necksize at most $2\pi/3$.  This bound in the symmetric case can now
be seen as the maximum side length for a spherical equilateral triangle.
More generally, we have the following~\cite{gbks1}.

\begin{corollary}
The triple $0 < x,y,z \le \pi$ can be the necksizes of a triunduloid if and
only if it satisfies the spherical triangle inequalities
$$ x+y+z\le 2\pi,\quad x\le y+z,\quad y\le z+x,\quad z\le x+y. $$
In particular, at most one end of a triunduloid can be asymptotic to a cylinder.
\end{corollary}

Similar methods apply to genus-zero surfaces with $k>3$ ends, when those
ends still have asymptotic axes in a common plane.  The moduli space of
such \emph{coplanar $k$-unduloids} can be understood as a covering of the
space of spherical $k$-gons.  More general surfaces, without coplanar ends,
will be more difficult to classify.


\end{document}